\documentclass[11pt]{amsart}\usepackage{amsfonts} \usepackage{amscd}
\usepackage{amssymb}

\newcommand{\ring}[1]{\mathbb{#1}}
\newcommand{\C}{\ring{C}} \newcommand{\Q}{\ring{Q}}
\newcommand{\Z}{\ring{Z}} 
\newcommand{\A}{\ring{A}}
\newcommand{\F}{\ring{F}} 
\newcommand{\N}{\ring{N}}

\newcommand{\lef}{\ring{L}}

\newcommand{\be}{\begin{equation}}
\newcommand{\ee}{\end{equation}}
\newcommand{\nd}{\noindent}

\newcommand{\arcs}{\mathfrak{L}}

\newcommand{\spec}{\mathnormal{\mathrm{Spec\,}}}
\def\1{{\mu\mkern-6mu\mu}}
\newcommand{\lie}{{\bf\mathfrak g}}

\newcommand{\G}{{\bf G}}

\newcommand{\so}{{\mathfrak {so}}}
\newcommand{\syp}{{\mathfrak {sp}}}

\newcommand{\all}{{\forall}}

\newcommand{\ord}{{\text{ord}}}

\def\rq{\text{'}}
\newcommand{\ri}{{\mathcal O}}

\author{Julia Gordon}
\address{Institute for Advanced Study, 1 Einstein Dr., Princeton
NJ 08540, USA}
\email{julygord@math.ias.edu}

\title[Motivic nature of character values]
{Motivic nature of character values of depth-zero representations}

\begin{document}

%\linespread{1.6}

\theoremstyle{plain}
\newtheorem{thm}{Theorem}
\newtheorem{lem}[thm]{Lemma}
\newtheorem{cor}[thm]{Corollary}
\newtheorem{prop}[thm]{Proposition}

\theoremstyle{definition}
\newtheorem{rem}[thm]{Remark}
\newtheorem{defn}[thm]{Definition}
\newtheorem{ex}[thm]{Example}

\footnotetext[1]{supported in part by the Fields Institute and NSF.}
%%%%%%%%%%%%%%%%%%%%%%%%%%%%%%%%%%%%%%%%%%%%%%%%%%%%%%%%%%%%%%%%%%
\newcommand{\ab}{\lambda}
\newcommand{\rf}{\F}
\newcommand{\val}{{\ord}}
\newcommand{\tf}{{\text{Tr\,Frob}}}
\newcommand{\mot}{\hat K_0^{mot}({\text {Var}}_F)_{\Q}}
\newcommand{\ac}{\overline{\text{ac}}}
\newcommand{\mx}{{\mathcal X}}
\newcommand{\rs}{{\text {Res}}}

\begin{abstract}
In the present paper, it is shown that the values of Harish-Chandra 
distribution characters on definable compact subsets of the set of topologically unipotent elements 
of some reductive $p$-adic groups can be expressed as the trace of Frobenius action on
certain geometric objects, namely, Chow motives.
The result is restricted to a class of depth-zero representations of symplectic or special orthogonal
groups that can be obtained by inflation from Deligne-Lusztig representations.
The proof works both in positive and zero characteristic, and relies on arithmetic motivic integration. 
\end{abstract}
\maketitle
\section{Introduction}

Our goal in this paper is to relate the 
distribution characters of depth-zero representations
of $p$-adic groups with
geometric objects, namely, Chow motives.
This is part of an effort, initiated by T.C. Hales in \cite{tomtalk}, 
to express the concepts of representation theory of $p$-adic groups
in such a way that would allow computations to be done 
without relying on the 
knowledge of the specific value of $p$ 
(see also \cite{TH1}, \cite{toi} and \cite{tf}).

Let $G$ be a $p$-adic group, for example, $G=\G(\Q_p)$ for a connected 
reductive algebraic group $\G$, and let $\pi$ be a representation of $G$.
Harish-Chandra \cite{HC} introduced the notion of
the {\it distribution character} $\Theta_{\pi}$ of $\pi$ and showed that
it is represented by a locally integrable function $\theta_{\pi}$ on the
group $G$. It turned out to be a difficult 
(and still unsolved in most cases) 
problem to give an explicit formula for the function $\theta_{\pi}$.
There is  evidence \cite{KLB}, \cite{non-elem} suggesting that
this difficulty is due to the existence of geometric objects that turn out
to be {\it non-rational}
whose number of points over the finite field $\F_p$ appears in the calculation
of the value of the character. 
This demonstrates that harmonic analysis on reductive groups is, in general,
{\it non-elementary} (see \cite{non-elem}), and the best we can 
hope for is not to obtain explicit formulas for character values, 
but  to understand the underlying geometry.
The first step towards this objective would be to 
establish the existence of geometric 
objects related to the values of characters in general.
This is the step we carry out here 
for a class of depth-zero representations.
The methods we use are completely different from the ones 
in \cite{non-elem}, and based on the new theory of {\it arithmetic motivic integration}.

In 1995, M. Kontsevich introduced the idea of motivic integration, which 
led J. Denef and 
F. Loeser to develop the theory of arithmetic motivic integration.
This theory  is outlined briefly in the Appendix.
Arithmetic motivic integration  
allows to express the classical $p$-adic volumes
of {\it definable} subsets of $p$-adic manifolds in terms of 
{\it trace of 
$p$-Frobenius action on virtual Chow motives} (which is, essentially, a 
generalization of counting $\F_p$-points of algebraic varieties).
In \cite{TH1}, T.C. Hales outlined a program that applies 
arithmetic motivic 
integration in order to relate various quantities arising 
in representation 
theory of $p$-adic groups (such as orbital integrals) to geometry.

Here we use these ideas to express some  averaged values of the 
 character
$\theta_{\pi}$ in terms of the Frobenius action on a Chow motive
associated with the group, the representation, and the set of averaging.
So far, such a result is restricted to the class of depth-zero 
supercuspidal representations of symplectic or special orthogonal (odd) 
groups, which can be obtained by inflation from Deligne-Lusztig 
representations. However, we expect that the methods should generalize 
to a much wider class, possibly including most supercuspidal 
representations.

{\bf Acknowledgment.} 
This work is part of a thesis written under the 
guidance of T.C. Hales, to whom I am deeply grateful for 
all the advice and help.
I am indebted to Ju-Lee Kim for 
suggesting the key idea of using the 
double cosets at the early stage of this work. Special thanks to 
J. Korman for all his help in getting this paper to its present shape, 
to J. Diwadkar for sharing her calculations with me, and to J. Adler and 
L. Spice for helpful conversations. I am also very grateful to the referee for
suggesting multiple improvements.

\section{The statement}\label{theorem}
\subsection{Notation and assumptions}\label{notation}
Throughout the paper, we assume that the algebraic group 
${\bf G}$ is a 
symplectic group $\G={\bf Sp}(2n)$ or special orthogonal group
$\G={\bf SO}(2n+1)$ for some $n\in\N$. 
Specifically, we think of ${\bf Sp}(2n)$, resp., ${\bf SO}(2n+1)$ as
the closed subgroup of ${\bf GL}(2n)$, resp., ${\bf GL}(2n+1)$,
cut out by the condition ${}^tgJg=J$.
 The matrix $J=(J_{ij})$ of the size $2n$ (resp., $2n+1$), 
in the symplectic case is defined by:
  $J_{ij}=0$ if $i+j\neq 2n+1$, $J_{i,2n+1-i}=1$ if $i\le n$,
and $-1$ otherwise, and in the orthogonal case  $J$ is the anti-diagonal
matrix with $1$'s on the anti-diagonal. The reason we have to restrict our 
attention only to these groups is explained in Section \ref{delu}.

The letter `$F$' will be reserved for a global field, and the letter
`$E$' -- for a local field. The symbol $\ri$ with various subscripts 
will always be used to denote the ring of integers of the corresponding field. 

We will  denote by 
${\mathcal K}$ the collection of all nonarchimedean completions
of a global field $F$.
Namely, we consider two cases: 
if $F$ is a number field, then  
${\mathcal K}=\{F_v\}$ is the collection
of its completions  at nonarchimedean places.
If $F$ is a function field, we denote by ${\mathcal K}$ the collection
of all fields of the form $\F_x((t))$, where $x\in\spec \ri_F$ is a
closed point.

If $k$ is a finite field, and $\G$ -- an algebraic group as above, we denote by 
$\G^k$ the finite group of $k$-points of $\G$, that is, the subgroup  of
$\G(\bar k)$ consisting of the points fixed by the Frobenius action (see, e.g.,
\cite[Chapter I]{Sri} for details).
  
\subsection{The representations}\label{inflation}
Let ${\mathcal K}$ be a collection of local fields, as defined in 
Section \ref{notation}.  
We would like to study the characters of representations of the groups
$\G(E)$, as $E$ ranges over the  family ${\mathcal K}$. 
As stated in the introduction, the goal is to obtain 
an independent of $E$ formula for the 
character. In particular, we need the 
representations to be, so to say, 
the ``incarnations of the same representation'' as we let $E$ range over
various completions of a given global field.
To make this precise, we consider a family of representations
of the groups $\G(E)$ (one representation for each group) that are 
parametrized by the same combinatorial data independent of $E$.

For the group $\G$, $\G={\bf SO}(2n+1)$ or $\G={\bf Sp}(2n)$ as above,
the starting datum is a partition $w$ of the integer $n$.
For each finite field $\F_{q}$, $w$ corresponds 
to a conjugacy class of elliptic tori in the finite group $\G^{\F_q}$
\cite[Section 5.7]{W2}.
In particular, 
for each field
$E\in{\mathcal K}$ with ring of integers $\ri_E$ and 
residue field $k_E$, 
$w$ corresponds to a conjugacy 
class of tori in $\G^{k_E}$.
We pick a torus  $T=T_{k_E,w}$ from this class. 
Deligne-Lusztig theory associates a representation $R_{T, \chi}^{\G^{k_E}}$ 
with
the data $(T, \chi)$ where $\chi$ is a character of the torus $T$.
We let the character $\chi$ be an arbitrary irreducible character of $T$
in general position.
Such character exists if the characteristic of $k_E$ is large enough; 
it will be explained in 
Section \ref{remarks} why the choice of $\chi$ subject to these restrictions 
does not matter. 
Let $\psi=R_{T, \chi}^{\G^{k_E}}$ be the 
Deligne-Lusztig representation of $\G^{k_E}$
corresponding to the pair $(T, \chi)$ 
(see Appendix, Section \ref{delu}). 
It is a representation of $\G^{k_E}$ on a $\bar \Q_l$-vector space
with $l\neq \text{char} k_E$, but its construction does not depend either 
on $k_E$ or on $l$. 
We would like to be able to vary
$E$ (and with it, $k_E$). 
Therefore, let us fix, once and for all, $l=5$ (see Section \ref{delu}
for the explanation of this choice),
and  only consider the local
fields in ${\mathcal K}$ 
with the residual characteristic greater than $5$ from now on.
We also embed $\bar \Q_5$ in $\C$, and replace $\psi$ with the
representation $\bar \psi$ on the resulting complex vector space.

The representation $\pi_{E,w, \chi}$ of the group $\G(E)$ associated with the 
partition $w$ and the character $\chi$ is constructed as follows.
Let $K_E=\G(\ri_E)$.
The group $K_E$ is a maximal compact subgroup of $\G(E)$.
First, we inflate the representation $\bar\psi$ to a representation
of $K_E$, and then induce the 
resulting representation from $K_E$ to $G(E)$
(by means of compact induction).
By `inflation' we mean the following process.
Denote by  
${\text {Res}:K_E\to\G^{k_E}}$ the map that acts, coordinate-wise, as 
reduction modulo the uniformizer of the field $E$. 
Then define the representation $\kappa:K_E\to \text{End}(V)$
by $\kappa(g):=\bar\psi(\text{Res}(g))$, where $V$ 
is the representation space of $\bar\psi$.
The final result of this inflation-induction procedure
is a depth-zero representation $\pi_{E,w, \chi}=\text{c-Ind}_{K_E}^{\G(E)}\kappa$ 
of the group
 $\G(E)$ on a complex 
vector space.

It is relevant to note here that the construction so far is not quite independent
of the field $E\in{\mathcal K}$, because of the presence of $\chi$.
However, later we will see that $\chi$ drops out of all calculations of the character
values. For now, it may be more precise 
to talk about a ``pack'' of representations $\pi_{E,w,\chi}$ associated with the
data $(E,w)$. 

The following known fact is important for us, so we state it as a lemma, 
but do not include the proof (cf. \cite[Proposition 6.8]{MP}).
\begin{lem}\label{sc}
The representation $\pi_{E,w, \chi}$ is supercuspidal.
\end{lem}

\subsection{The main theorem}\label{statement}
Let $\G$, ${\mathcal K}$, $w$, $\chi$, and $\pi_{E,w, \chi}$ for almost all 
$E\in {\mathcal K}$ be as in Section \ref{inflation}.
For almost all $E\in {\mathcal K}$ the  
residual characteristic of $E$ is large enough so that the construction 
of the Section \ref{inflation} makes sense.

We denote the distribution character of the representation $\pi=\pi_{E,w, \chi}$
by $\Theta_{w,E}$, even though it would be more precise to 
denote it by $\Theta_{w,\chi,E}$. We are omitting $\chi$ from the notation 
because our calculations do not depend on $\chi$.  
We denote  the locally summable function that represents the 
distribution $\Theta_{w,E}$ by $\theta_{w,E}$, 
and also call it the character of $\pi$. 

Let $K_E^{rtu}$ denote the set of regular  topologically unipotent
elements in $K_E=\G(\ri_E)$.
Denote by $I_E$ the Iwahori subgroup of the group $\G(E)$ which 
consists of the elements in $K_E$ whose reduction 
modulo the uniformizer is an upper-triangular matrix in the 
standard representation.

To state the main theorem, we need two notions: of Pas's language, and
of virtual Chow motives. They are defined in the Appendix: 
Sections \ref{pl}, \ref{def}, and \ref{mots}.
The completed ring of virtual Chow motives which ``arise from varieties'' 
(see Apenndix, Section \ref{mots}) is denoted by $\mot$.
 
Let $\alpha$ be a formula in Pas's language
(see Appendix, Section \ref{def}) defining 
a compact subset $\Gamma_{\alpha,E}=Z(\alpha,E)\subset K_E^{rtu}$ 
for almost all 
$E\in{\mathcal K}$.
Let $vol$ stand for 
the $p$-adic Haar measure on $G=\G(E)$ normalized so that its
restriction
to $K_E$ coincides with the Serre-Oesterl\'e measure 
(see Appendix, Section \ref{armint}). 

\begin{thm}\label{weak}
Given $\G$, $w$ and ${\mathcal K}$ as above,  
there exists a virtual Chow
motive $M_{\alpha,w}\in\mot$, such that for 
almost all $E\in {\mathcal K}$,
the following equality holds:
$$
 vol(I_E) \Theta_{w,E}(\Gamma_{\alpha,E})
=vol(I_E)\int_{\Gamma_{\alpha,E}}\theta_{w,E}(\gamma)\,d\gamma
  =TrFrob_E(M_{\alpha,w}),
$$
where $Frob_E$ is the Frobenius action corresponding to the place of $F$ 
that gave rise to the  field $E$, and
$\Theta_{w,E}(\Gamma_{\alpha,E})$ stands for the value of 
the distribution $\Theta_{w,E}$ at the characteristic function of the set
$\Gamma_{\alpha,E}$.
\end{thm}

As announced in the introduction, this theorem states that  
the value of the distribution character on a {\it definable} 
compact subset of $K^{rtu}$ can be recovered from a geometric 
object (namely, a virtual Chow motive), up to a factor 
which is a polynomial in the 
cardinality of the residue field (the 
volume of the Iwahori subgroup).

\subsection{Remarks}\label{remarks}
{\bf 1.} A new version of the theory of motivic integration, that is 
about to appear, yields that the ring of virtual Chow motives does 
not need to be completed in order to define the elements $M_{\alpha, w}$
(see Appendix, Section \ref{mots}). 

{\bf 2.}
The reason we restrict our attention only to the topologically 
unipotent elements is the following.
Recall that the irreducible character $\chi$ of the torus in the 
finite group $\G^{k_E}$ is part of the construction of the 
representation $\pi$ which does not have a nice combinatorial 
parametrization.
However, as we will see in the next section, the 
value of 
$\theta_{w, E}$ at a {\it topologically unipotent} element 
is expressed through the values of the character of Deligne-Lusztig 
representation $\bar\psi$ that gave rise to $\pi$ at {\it unipotent}
 conjugacy classes. Those values, in turn, are expressed by Green polynomials
(see Appendix, Section \ref{delu}, for a brief review and references). 
We are using the fact that 
the values of the Green polynomials at 
unipotent elements
do not depend on the choice of the character $\chi$.

{\bf 3.}
We are considering  the averages of the function 
$\theta_{w,E}$ over some definable subsets of
$K_E^{rtu}$.
The other way of saying this is to say that we look at the values
of the distribution character  $\Theta_{w,E}$ itself 
at the characteristic functions of those subsets.
 In fact, we would like to study
the individual values  $\theta_{w,E}(\gamma)$ at regular
elements. Pas's language, which is one of our  main tools, does not allow
to handle individual elements of $p$-adic fields or  $p$-adic groups. 
As a  function
on the set of regular elements in the $p$-adic group $\G(E)$, the 
character 
$\theta_{w,E}(\gamma)$ is locally constant.
However, there are no explicit results  
that say ``how small'' the sets
on which $\theta_{w,E}$ is constant might be. 
This forces us to average
$\theta_{w,E}$ over some sets that we can control. 
However, the actual shape of
the sets  of averaging is flexible 
(the only requirement being that
they are {\it definable}, see Section \ref{def}).

The rest of the paper is devoted to the proof of Theorem \ref{weak}.
In the next section we do an entirely $p$-adic calculation which 
serves as a preparation for the proof, and in Section \ref{moti}
we give it a motivic interpretation, which leads to the desired 
result.

\section{The Character: a $p$-adic calculation}
Throughout this section we adhere to a fixed local field $E$ with
a uniformizer $\varpi$, ring of integers $\ri$, and   residue 
field $k=\ri/(\varpi)$.
Since the field $E$ stays fixed here, we drop the subscript $E$
from all notation for simplicity.

\subsection{Character of an induced representation}\label{charind}
Let $\pi$ be a representation of $G=\G(E)$ 
obtained from a Deligne-Lusztig
representation $\bar\psi$ of $\G^k$ by the 
inflation and induction  procedure described
in  Section \ref{inflation}. 
By Lemma 1, $\pi$ is supercuspidal.

We fix $K:=K_E$ -- a maximal compact subgroup of $G$, and 
let $I=I_E$.

We denote by $\rho$ the character of the representation $\bar\psi$, 
and by  
 $\tilde\rho\colon K\to\C$ -- the character
of its inflation to $K$.
Note that the value  $\tilde\rho(g)$ at
$g\in K$
depends only on the reduction of $g$ modulo $\varpi$.

Let $\Gamma$ be a compact subset of $K^{rtu}$.
By the Frobenius-type formula for induced character, 
\cite [Theorem 1.9]{sally} and \cite[Theorem  A.14(ii)]{bh} (see also the Remark at the end of the Appendix in \cite{bh}),
the value of the character $\theta_w(\gamma)$
can be expressed as a sum that is {\it finite} uniformly for all 
$\gamma\in\Gamma$: 
\begin{equation}\label{induced}
\theta_w(\gamma)=
\sum_{a\in I\backslash G/K}\left(\sum_{y\in IaK/K}
\tilde\rho(y^{-1}\gamma y)\right).
\end{equation}
Note that in \cite{sally}, 
the double cosets are taken with respect to the same 
subgroup on the left and right; however, the modification we are using here
can be obtained in the exactly same way. Our version of the 
formula coincides with the formula in \cite{bh}. 

\subsection{The double cosets}\label{dc}
The indexing set $I\backslash G/K$ in 
the outer sum in (\ref{induced}) 
can be described using a version of the Iwasawa
decomposition.

Let $A$ be the split standard torus in $G$. We can think of the 
elements of $A$
as diagonal matrices of the form
$\text{diag}(u_i\varpi^{\lambda_i})_{i=1,\dots,r}$
with $u_i\in \ri^{\ast}$ and $\lambda_i\in\Z$ satisfying the relations
forced by the definition of the group (for example, in the symplectic
case, $u_i=u_{r-i}^{-1}$ and $\lambda_i=-\lambda_{r-i}$). Then, by a
version of Iwasawa decomposition,  $G=IAK$ \cite{I1}.

Observe that
the intersection of the subgroups $I$ and $A$ is exactly the
intersection of $A$ with $K$. It follows that every left
coset of $G$ modulo $K$ has a representative of the form $ya$ with $y\in
I$ and $a\in A_0$, where $A_0$ is the set of
 elements in $A$ of the form
$a_{\ab}=\text{diag}(\varpi^{\ab_i})_{i=1,\dots,r}$, where
$\ab=(\ab_i)_{i=1,\dots,r}\in\Z^r$ satisfies the conditions mentioned above.
In particular, the set of double cosets $I\backslash G/K$ is in 
bijection with the set $A_0$.

Fix a multi-index $(\lambda_i)_{i=1,\dots,r}$ that gives an element
$a_{\lambda}\in A_0$.  
\begin{defn}
We call the elements $y_1, y_2\in I$ $\lambda$-equivalent
if $y_1 a_{\lambda}$ and $y_2 a_{\lambda}$ belong to the 
same left $K$-coset
(that is, $a_{\ab}^{-1}y_1^{-1}y_2 a_{\ab}\in K$).
We denote the $\lambda$-equivalence class of $y\in I$ by $[y]_{\lambda}$.
\end{defn}

As an element of the torus $A$, $a_{\lambda}$ may be considered 
as a lift of an element
of the extended Weyl group of ${G}$. Let $l_{\lambda}$ be the length 
of this element. 
By a theorem of Iwahori and Matsumoto \cite{I1},
the set of $\ab$-equivalence classes is in bijection with
$\A^{l_{\ab}}(k)$.

The following two simple observations will be used below, 
so we state them as a lemma. We keep the same notation as above, 
and let $q=|{\mathcal O}/(\varpi)|$ denote the cardinality of the 
residue field.

\begin{lem}\label{eqcl}
{\bf 1.} If $y_1$ and $y_2$ are $\ab$-equivalent, then
the elements \newline
$(y_1 a_{\ab})^{-1}\gamma (y_1 a_{\ab})$ and
$(y_2 a_{\ab})^{-1}\gamma (y_2 a_{\ab})$ are 
conjugate by an element of $K$.
\item{}{\bf 2.}
All the $\ab$-equivalence classes have equal volumes, equal to 
$\frac {vol(I)}{q^{l_{\ab}}}$.
\end{lem}\nd
{\bf Proof.} The first statement is obvious:
$
(y_1 a_{\ab})^{-1}\gamma (y_1 a_{\ab})$ \newline
$=k (y_2 a_{\ab})^{-1}\gamma (y_2a_{\ab}) k^{-1}$,
where $k=a_{\ab}^{-1}y_1^{-1}y_2 a_{\ab}\in K$ by definition of 
$\ab$-equivalence.

To show the second assertion, note that 
$y$ and $y_0$ are $\ab$-equivalent if and only if 
$y^{-1}\in a_{\ab}K a_{\ab}^{-1}y_0^{-1}\cap I$. 
Hence the volume of each equivalence class equals 
the volume of the set $I\cap a_{\ab}Ka_{\ab}^{-1}$ (since the group $G$ is 
reductive, it is unimodular, and therefore
the operation of taking an inverse preserves the Haar measure).
Let us compute this volume. The
total number of equivalence classes within the Iwahori subgroup $I$
equals  $q^{l_{\ab}}$. Hence, the volume of each equivalence class is
$\frac {vol(I)}{q^{l_{\ab}}}$.
\qed

\subsection{A formula for the character}\label{f}
First, let us introduce some notation. 
Let $C$ be a conjugacy class in the finite group $\G^k$, and let 
$\gamma\in K^{rtu}$.
For $\lambda\in\Z^r$ as in the previous section,
denote by $N_C^{\ab}(\gamma)$ the number of $\ab$-equivalence classes
$[y]_{\ab}$
of elements of $I$ such that the following conditions are satisfied:
\newcounter{cond}
\begin{list}
{\roman{cond})}{\usecounter{cond}}
\item\label{fix}
$(ya)^{-1}\gamma (ya)\in K$  for any $y\in [y]_{\ab}$
\item\label{conj}
$\rs((ya)^{-1}\gamma(ya))\in  C$ for any $y\in [y]_{\ab}$.
\end{list}
The condition (i) is well defined by Lemma \ref{eqcl}.
Let $N_C(\gamma)=\sum_{\ab}N_C^{\ab}(\gamma)$.
We observe that for each element $\gamma\in K^{rtu}$ the sum is finite, 
because the set of classes $[y]_{\ab}$ satisfying the condition
(i) is nonempty only for a finite set of multi-indices $\ab$,
by \cite[Theorem A.14(ii)]{bh}.

Using  the 
fact that $\tilde \rho$ is a lift of
the character $\rho$ of the finite group $\G^k$,
 the sum (\ref{induced}) can be rewritten as
a sum over the  conjugacy classes $C$ of $\G^{k}$:
\begin{equation}\label{sum}
\theta_w(\gamma)=\sum_{C}\sum_{(\lambda, [y]_{\lambda}) \atop{{\all y\in
[y]_{\ab} (ya)^{-1}\gamma(ya)\in K,} \atop{\rs((ya)^{-1}\gamma(ya))}\in
C}}\rho(C)=\sum_C N_C(\gamma)\rho(C),
\end{equation}
where $C$ runs over the conjugacy classes of $\G^{k}$.

\begin{lem}\label{fin}
{\bf 1.} If $\Gamma$ is a compact subset of the set of regular elements
in $G$, then there exists a finite set $A_{\Gamma}$ of indices $\lambda$
so that $N_C^{\ab}(\gamma)=0$ for all $\gamma\in\Gamma$
and $\ab\notin A_{\Gamma}$.
\item{}{\bf 2.} 
If $\gamma$ is topologically unipotent, only the unipotent conjugacy classes
$C$ appear in the summation in (\ref{induced}). 
\end{lem}
\nd
{\bf Proof.} The first statement is part of the statement of 
\cite[Theorem A.14(ii)]{bh}. 
(In \cite{bh}, the Theorem is stated only for $GL_n$ and its relatives, but
as the authors point out at the end of the Appendix, the argument
carries over to classical groups without any changes.)

The second statement  is based on the following two simple observations.
First, the eigenvalues of $\text{Res}(\gamma)$ are obtained from 
the eigenvalues of $\gamma$ by reduction modulo 
the uniformizer of the valuation extended to the algebraic closure of
the field $E$. 
Second, it follows, for example, from \cite[Lemma 2.2]{Cassels}
that the eigenvalues of a topologically unipotent element
(in the algebraic closure of the
given local field) are congruent to $1$.
\qed
 
\subsection{Averages}\label{av}
Let $\Gamma$ be a compact subset of $K=\G(\ri)$.
The goal
is to express the average $\Theta_{w}({\Gamma})$ of the values
of $\theta_{w}$ over the set $\Gamma$ as the $p$-adic  volume of some
$p$-adic object. In the next section we will use Denef and Loeser's
``comparison theorem'' to relate the resulting $p$-adic volume to a Chow
motive.
We do it by means of artificially constructing some subsets of
$G\times G$ whose $p$-adic volumes equal the numbers $N_C(\gamma)$ that
appear in  formula (\ref{sum}), and showing that these subsets are 
definable.

Consider the Cartesian product
$$
\mx=\G\times_{\ri} \G.
$$
Fix an element $a=a_{\ab}$ as above and consider the double 
coset $D_a=IaK\subset G$.
Let $W_{\ab}^E(\Gamma)$ be the subset of $\mx(\ri)$
defined by
$$
W_{\ab}^E(\Gamma)=\{(y,\gamma)\mid y\in I, \gamma\in {\Gamma},
a^{-1}y^{-1}\gamma ya\in K\}.
$$

The set $W_{\ab}^E(\Gamma)$ is partitioned into subsets
$$
W_{C,\ab}^E(\Gamma)=\{(y,\gamma)\in W_{\ab}^E(\Gamma)\mid 
\rs({(ya)^{-1}\gamma (ya)})\in C\},
$$
where $C$ runs over the unipotent conjugacy classes of $\G^{k}$.

Let $W_{C,\ab}^E(\gamma)\subset I$ be the projection onto the first 
coordinate of the
cross-section  of $W_{C,\ab}^E(\Gamma)$ with fixed $\gamma\in\Gamma$.

Let $vol$  be the Haar measure on $\G(E)$ normalized so that its
restriction
to $K$ coincides with the Serre-Oesterl\'e measure. 
Then the  space $\mx(\ri)$ has the
natural product measure.

\nd{\bf{Main Observation.}} Up to a normalization factor,
the number $N_C^{\ab}(\gamma)$ defined in Section \ref{f}
is the volume of the set $W_{C,\ab}^E(\gamma)$:
$$
N_C^{\ab}(\gamma)=\frac{q^{l_{\ab}}}{vol(I)}vol(W_{C,\ab}^E(\gamma)),
$$
where 
$q$ is the cardinality of the residue field $k$, $l_{\ab}$ is an 
integer that depends only on $\lambda$, and $I$ is the Iwahori subgroup.

\nd{\bf Proof.} Fix a multi-index $\ab$. By definition,
$N_C^{\ab}(\gamma)$ is the number of equivalence classes
$[y]_{\ab}$ of the elements $y\in I$ that satisfy the conditions (i) and (ii).
The set $W_{C, \ab}^E(\gamma)$ is a disjoint union of these equivalence
classes. 
Recall that by Lemma \ref{dc}, all $\ab$-equivalence classes have equal 
volumes, equal to $\frac {vol(I)}{q^l_{\ab}}$.
Hence,
$$
N_C^{\ab}(\gamma)=\frac{q^{l_{\ab}}}{vol (I)}
vol(W_{C, \ab}^E(\gamma)).
$$
In the next section, we will also need the following property of the 
normalization factor that appears in the above formula.

\begin{rem}
The factor  $p(q)=\frac{q^{l_{\ab}}}{vol(I)}$, as a  function of
$q$, is an element of $\Z(x)$,  since the volume
of $I$ is a polynomial in $q$ with coefficients in $\Z$ 
that depends only on the group $\G$.
\end{rem}

Consider the average (with respect to the Haar measure 
$vol$ defined above)
of the values of $\theta_{w}$ over the set
$\Gamma$. In order to get rid of denominators, we multiply it by 
$vol(I)$. 
\begin{align}\label{aver}
vol(I)\Theta_{w}({\Gamma})&=
vol(I)\int_{\Gamma}\theta_{w}(\gamma)\,d\gamma=
vol(I)\int_{\Gamma}\sum_{C}\sum_{\ab}N_C^{\ab}(\gamma)\rho(C)\\
\label{aver1}&=
\sum_C\sum_{\ab}\rho(C)q^{l_{\ab}}\int_{\Gamma}vol(W_{C,\ab}^E(\gamma))
\,d\gamma\\
\label{aver2}&=
\sum_C\sum_{\ab}\rho(C)
q^{l_{\lambda}}\iint_{W_{C,\ab}^E(\Gamma)} 1\, dy d\gamma=
\sum_C\sum_{\ab}\rho(C)q^{l_{\ab}}
vol({W}_{C,\ab}^E(\Gamma)).
\end{align}
Note that the sum and the integral in the expression (\ref{aver1})
can be interchanged because the summation over $\lambda$, in fact,  
runs over  a finite set 
$A_{\Gamma}$, by Lemma \ref{fin}.  
This equality is the basis for our main result, as the next
section shows.

\section{Motivic interpretation}\label{moti}
\subsection{Varying the field}\label{vf} 
In this section, we complete the proof of Theorem \ref{weak} by 
associating virtual Chow
motives with the subsets $W_{C,\ab}^E(\Gamma)$ from the previous section
using arithmetic motivic integration. This
leads to an expression of the average value of $\theta_{w}$ 
as a trace of
Frobenius on a certain virtual Chow motive.

Throughout this section, we work with the notation and all the assumptions 
of Theorem \ref{weak}, so that ${\mathcal K}$ is the collection of all 
nonarchimedean  completions
of a given global field $F$; $\G$, $\pi$, $w$ are as in 
Section \ref{inflation}, 
$\alpha$ is a formula in Pas's language, and 
for $E\in{\mathcal K}$, $\Gamma_{\alpha,E}$ is a
definable subset of $K_E$ defined by the formula $\alpha$.
As in Theorem \ref{weak}, we are assuming that 
$\Gamma_{\alpha,E}$ is contained in the set of regular topologically
unipotent elements in $K_E$. 
In order to justify the assumption that a {\it definable} set 
$\Gamma_{\alpha,E}$
can be  contained in $K^{rtu}$, we would like  
to show that the set $K^{rtu}$ itself is 
 definable.

\begin{lem}\label{alldef}
Let $E\in{\mathcal K}$ be a local field. Let $K=K_E$.
Then
\begin{enumerate}
%\counter{cond}\setcounter{1}
\item[1.]{}
The set $K^{rtu}$ is definable.
\item[2.]{}
The Iwahori subgroup $I_E$ is a definable subset of $K_E$.
\end{enumerate}
\end{lem}\nd
{\bf Proof. 1.} Let $\gamma$ be an element of $K^{rtu}$.
We use the standard representation of the group $G$ to think of
$\gamma$ as a matrix $\gamma=(x_{ij})$, and let the `$x_{ij}$' be the 
free variables in all our logical formulas.
We use the following definition of {\it regular} \cite{KN}.
Let $n$ be the dimension of $\G$ and $l$ be the rank of $\G$.
Then the element $\gamma$ is regular if and only if
$D_l(\gamma)\neq 0$, where $D_l(\gamma)$ is defined by
the following expression: 
$$\det((t+1)I_n-\text{Ad}(\gamma))
=t^n+\sum_{j=0}^{n-1}D_j(\gamma)t^j.$$
We observe that $\text{Ad}(\gamma)$ can be thought of as
 a matrix expression whose entries 
are polynomials in $x_{ij}$, and hence $D_l(\gamma)$ is also a polynomial
expression in the variables $x_{ij}$.
Therefore, `$D_l(\gamma)\neq 0$' is a formula in Pas's language.
In conjunction with the formulas coming from the equations of the
group $\G$ as an affine variety, and with the Pas's language formulas
`$\val(x_{ij})\ge 0$' for each index $(i,j)$, this formula defines
the set of regular  elements in $K$.

To cut out the subset of topologically unipotent elements,
we use Lemma \ref{fin}, which implies that $\gamma$ is topologically unipotent
if and only if the element $\rs(\gamma)$ is unipotent, 
and then use the combinatorial (independent of the field) parametrization
of the set of unipotent conjugacy classes in the finite group $\G^k$
(where $k$ is the residue field of $E$).
The details are part of the proof of Lemma \ref{conjcl} below.

{\bf 2}. The Iwahori subgroup $I$ that we are considering is defined by
the formulas `$\val(x_{ij})\ge 0$', $i\le j$, 
in conjunction with the formulas
`$\val(x_{ij})\ge 1$' for all pairs of indices $(i,j)$ with $i>j$.
\qed

%%%%%%%%%%%%%%%%%%%%%%%%%%%%%%%%%%%%%%%%%%%%%%%%

\subsection {Preparation: the subsets $W_{C,\ab}^E(\Gamma_{\alpha,E})$ are definable}\label{prep}
We start by fixing a multi-index $\ab$, and considering
 the corresponding double coset $D_a$, as
in Section \ref{av}. Let $l=l_{\ab}$.
Recall that $\alpha$ stands for a formula in Pas's language that defines a compact subset
$\Gamma_{\alpha, E}$ of $K_E^{rtu}$ for almost all $E\in K$. Since we are keeping $\alpha$ 
fixed, in this section we drop `$\Gamma_{\alpha, E}$' from the notation, and denote the sets
whose definability we wish to prove simply by 
$W_{C,\ab}^E$. 

\begin{lem}\label{integr}
Let $\ab$ be a fixed index.
Under the assumptions of Theorem \ref{weak}, 
$W_{\ab}^{E}(\Gamma_{\alpha,E})$ is a definable subset of 
$\G\times \G(E)$ for $E\in{\mathcal K}$. 
\end{lem}\nd
{\bf Proof.}
Suppose for a moment that the field $E\in{\mathcal K}$ is fixed,
an element $y=(y_{ij})\in I_E$ is fixed, and
take $\gamma=(\gamma_{ij})\in\Gamma_{\alpha,E}$.
We observe that each matrix coefficient of the matrix
$y^{-1}\gamma y$ is a polynomial, with coefficients in $\Z$, homogeneous
in the  variables $y_{ij}$ and linear in the variables $\gamma_{ij}$.
We denote these polynomials by $P_{\kappa\eta}(y_{ij},\gamma_{ij})$,
$\kappa,\eta=1,\dots,r$.
The conjugation by $a_{\ab}$ multiplies each entry of $y^{-1}\gamma y$
by an integral power of the uniformizer
that depends only on the multi-index $\lambda$.
We denote these powers by $n_{\kappa\eta}$.

Consider the formulas
\begin{equation}\label{phis}
`\val P_{\kappa\eta}(y_{ij},\gamma_{ij})+
n_{\kappa\eta}\ge 0\rq,\quad `\val(y_{ij})\ge 0\rq, i\le j, \quad
`\val(y_{ij})\ge 1\rq, i>j.
\end{equation}

Let $\phi$ be the conjunction of all the  formulas (\ref{phis}) and the
formula  $\alpha$ in the variables $\gamma_{ij}$ that defines 
$\Gamma_{\alpha,E}$.
Then $W_{\ab}^{E}(\Gamma_{\alpha,E})=Z(\phi,E)$ in the notation of Section \ref{def} of the 
Appendix.
\qed

\begin{lem}\label{conjcl}
The subsets $W_{C,\ab}^{E}$ are definable.
\end{lem}\nd
{\bf Proof.}
Suppose for a moment that a local field $E\in {\mathcal K}$ is fixed.
 All we need to show is that
the set of topologically unipotent elements  $y\in K$ such that
 $\rs(y)\in C$ is definable for each conjugacy class
$C$ of the group $\G^k$ that appears in the summation (\ref{aver}).
In order to write a formula that cuts out this set, 
we let the
matrix coefficients of $y$ be the free variables, as before.

First of all, we observe that the symbol `$\rs$' can be used in formulas
in Pas's language with the following meaning:
if $\varphi(y)$ is an expression in Pas's language with $y$ -- a variable of the residue 
field sort, then we set $\varphi(\rs(x))$
(where $x$ is a free variable of the valued field sort) 
to be $\varphi_1\vee\varphi_2$, where
$\varphi_1:=$`$\val(x)=0$'$\wedge \varphi(\ac(x))$, and $\varphi_2:=$`$\val(x)>0$'$\wedge\varphi(0)$.
Also, if $x$ is an abbreviation for a matrix or a vector of free variables, we understand
`$\rs(x)$' in a natural way, as the application of the symbol `$\rs$' to each component 
of $x$, as described above.

Since $\Gamma_{\alpha,E}$, by assumption, is contained in the set of
topologically unipotent elements in $G$, only the {\it unipotent} 
conjugacy classes
$C$ appear in (\ref{aver}), by Lemma \ref{fin}.
As described in the Appendix, Section \ref{unip}, 
the set of unipotent conjugacy classes
$C$ is in bijection with the set of purely combinatorial data
$\{(\Lambda, (\epsilon_i))\}$, where $\Lambda$ is a partition of 
$r$, and $\epsilon_i\in k_E^{\ast}/{k_E^{\ast}}^2$ 
(where $k_E$ is the residue field of $E$).
Suppose that the class $C$ corresponds to a pair $(\Lambda, (\epsilon_i))$.
We think of $(\epsilon _i)$ as a sequence with entries $0$ or $1$.

In order to write down the logical formula $\phi_{\Lambda, (\epsilon_i)}$
cutting out the set of elements whose reductions fall into the conjugacy 
class $C$, we need to unwind the process described in  Section \ref{unip}
that associates $C$ with the data $(\Lambda, (\epsilon_i))$. 

First, for $y\in K$, $y$ -- topologically unipotent, 
let $Y=(1-y)(1+y)^{-1}$ be
the Cayley transform of $y$ (note that the matrix $1+y$ is invertible,
since $y$ is assumed to be topologically unipotent). The element $Y$ lies in 
the Lie algebra of $\G$ over the given $p$-adic field. 
Let $\tilde Y_1=\det(1+y)Y$, and let `$Y_1$'=`$\rs(\tilde Y_1)$'
(so that `$Y_1$' in all formulas is treated as a matrix with components ranging 
over the residue field sort, even though, to be precise, all the formulas involving `$Y_1$'
would be  conjunctions of formulas with free variables ranging over the valued field sort, namely,
the matrix coefficients of $\tilde Y_1$).
Note that the matrix coefficients
of $\tilde Y_1$ are polynomials, with
$\Z$-coefficients, in the matrix coefficients of $y$. 
Hence, finally,  the abbreviation `$Y_1$' can be used 
instead of `$y$' in all subsequent formulas, and we treat it just as a matrix 
with components of the residue filed sort. 
 
The first part of the data, $\Lambda$, prescribes the set of Jordan 
blocks of the matrix $Y$, which is the same as the set of Jordan blocks
of $Y_1$.
The set of all elements $Y_1$ whose Jordan blocks match the partition
$\Lambda$ is defined by
\begin{equation}\label{phi1}
`\exists
(g_{ij})_{i,j=1,\dots, r} : (g_{ij})Y_1(g_{ij})^{-1}=J_{\Lambda}\rq,
\end{equation}
where $J_{\Lambda}$ is the matrix with entries $0$ and $1$, consisting of 
Jordan blocks given by the partition $\Lambda$.
Note, again, that in fact formula (\ref{phi1}) has only free variables of the valued field sort, that
is, the matrix coefficients `$y_{ij}$' (we see that by unwinding the abbreviation `$Y_1$'). 
All the variables ranging over the residue field sort (that is, `$g_{ij}$') are bound.
Hence, 
logical formula (\ref{phi1}) cuts out a definable set.

Second, we need to cut out the subset of this definable set that
corresponds to the given sequence $(\epsilon_i)$. 
Recall (from Section \ref{unip}) that the
sequence $(\epsilon_i)$ comes from a collection of quadratic forms on
the vector spaces $V_i$. Let $c_i=\dim V_i$.
We will need to
introduce separately a few components of the final formula that would
define $W_{C, \ab}^E$. 

Let us look at one entry $\epsilon_i$ corresponding
to the vector space $(V_i, q_i)$.

Let `$\text{li}(v_1,\dots, v_m)$' be the logical formula with the free
variables $v_{11}, \dots, v_{mr}$ that states that the vectors
$v_1={}^t(v_{11}, \dots, v_{1r}), \dots, v_m={}^t(v_{m1}, \dots,
v_{mr})$ are linearly independent (here $r=\dim V$), see \cite{tf}.

Let `$\ker_A(v)$' be the logical formula stating that the vector $v$ is
annihilated by the linear operator $A$, 
where `$A$' is a matrix of logical terms. 
The free variables of this formula are the components of $v$ and 
the matrix coefficients of $A$.
Let `$\text{q-span}(Y_1,v_1,\dots,
v_m)$' be the formula with free variables 
$v_1, \dots, v_m$ and $Y_1$ (each of which is an 
abbreviation for a vector or a matrix, respectively), which is the conjunction
of the formulas $\ker_{Y_1^i}(v_j)$,
$j=1,\dots,m$ and the following:
\begin{align*}
`\all u \exists (a_1,\dots, a_m) \quad\wedge\exists v'\wedge \exists v''\\
\ker_{Y_1^{i-1}}(v')\wedge \ker_{Y_1^{i+1}}(v'')\quad\wedge
v'+Y_1v''+\sum_{i=1}^m a_i v_i = u\rq,
\end{align*}
where $a_j$ are scalars, and $v', v'', u$ are
abbreviations for  $r$-vectors (all existential quantifiers in this formula 
range over the residue field sort). This formula states that the vectors
$v_1, \dots, v_m$ span $\ker (Y_1^i)$ modulo $\ker
(Y_1^{i-1})+Y_1\ker(Y_1^{i+1})$. Let $J$ be the matrix of the quadratic form
$q_V$, that is, the matrix from the definition of $\G$ as a subgroup of 
${\bf GL}(r)$. 
Let
`$\text{Gram}(v_1,\dots, v_c)$' stand for the matrix of formal expressions
that gives the  Gram
matrix of the
basis $v_1,\dots, v_c$ with respect to the quadratic form 
defined by the matrix ${}^t(Y_1^{i-1}w)Jw'$.

If the
$\epsilon_i$ given is $1$, let `$\psi_i(y)$' be the formula:
\begin{align*}
\lq \exists
(v_1,\dots,v_{c_i})\ \wedge \ \text{q-span}(Y_1, v_1,\dots,v_{c_i})\\
\wedge \text{li}(v_1,\dots,v_{c_i})\ \wedge \
 \exists (z\in k\quad
z^2=\det(\text{Gram}(v_1,\dots, v_{c_i}))\ \wedge \ z\neq 0).
\rq
\end{align*}

If the
$\epsilon_i$ given is $0$, let `$\psi_i(y)$' be the formula:
\begin{align*}
\lq
\exists
(v_1,\dots,v_{c_i})\ \wedge \ \text{q-span}(Y_1, v_1,\dots,v_{c_i})\\
\wedge \text{li}(v_1,\dots,v_{c_i})\ \wedge \
 \nexists (z\in k\quad
z^2=\det(\text{Gram}(v_1,\dots, v_{c_i}))\ \wedge\ z\neq 0).
\rq
\end{align*}

Finally, it follows from Section \ref {unip} that 
the set $W_{C, \ab}^{E}$ is defined by the conjunction of $\psi_i(y)$
for all $i$ and the formula defining the set $W_{\ab}^E(\Gamma_{\alpha,E})$.  \qed

\subsection{Proof of Theorem \ref{weak}}\label{pff}
Let ${\mathcal K}$ be the collection of local fields as in the statement 
of the theorem. Let $E\in {\mathcal K}$.
By the formulas (\ref{aver}) -- (\ref{aver2}),
the value $vol(I_E)\Theta_w({\Gamma_{\alpha,E}})$
is a sum over $\ab$ of
the expressions of the form\newline
$\sum_C\rho(C)q^{l_{\ab}}
vol(W_{C,\ab}^{E}(\Gamma_{\alpha,E}))$. We are now ready to write down the corresponding 
motivic expression for each of these terms.

Let us fix $\ab$ for now.
Let ${\mathcal X}=\G\times \G$  as before.
By Lemmas \ref{integr}  and \ref{conjcl},
the sets ${{W}}_{C,\ab}^{E}(\Gamma_{\alpha,E})$ are definable subsets of $\mx(E)$, 
i.e., there exist 
formulas $\phi_{C,\ab}$ in Pas's language, such that
 $W_{C, \ab}^{E}(\Gamma_{\alpha,E})=Z(\phi_{C,\ab},E)$ (in the notation of Appendix, Section 
\ref{def}).
Let $M_{C,\lambda}=\mu(\phi_{C, \ab})\in\mot$
be the arithmetic motivic volume of the formula $\phi_{C,\ab}$
(see Appendix, Section \ref{armint}).

By \cite{Sr3}, $\rho(C)$ is a polynomial in $q$ for each $C$, where
$q$ is the cardinality of the residue field of $E$.
Denote this polynomial by $f_C(q)$, and let
$\tilde M_{C,\ab}=f_C(\lef)\lef^{l_{\ab}}M_{C,\ab}$, where 
$\lef$ is the Lefschetz motive (see Appendix, Section \ref{mots}), 
so that $\tilde M_{C,\ab}$ 
is also an element of the ring $\mot$.

Finally, let
$$
M_{\ab}=\sum_{(\Lambda, (\epsilon_i))}\tilde M_{C,\ab}\in\mot,
$$
where $(\Lambda, (\epsilon_i))$ is the data parametrizing the 
unipotent conjugacy classes
$C$, as in Section \ref{unip}.

By Theorem 8.3.1, \cite{DLar}
(respectively, Theorem 8.3.2 in the function field case, see also Appendix, 
Section
\ref{armint}),
\begin{equation}\label{l-dep}
q^{l_{\ab}}\rho(C)\iint_{W_{C,\ab}^E(\Gamma_{\alpha,E})} 1\, dy d\gamma=
 TrFrob_E(M_{C,\ab})
\end{equation}
for almost all completions $E$ of $F$, and hence the sum of the expressions
on the left-hand side of (\ref{l-dep}) over $C$ equals $TrFrob_E(M_{\ab})$ 
for almost all $E$.

In order to obtain the final result, it remains to sum up both sides
of the above formula over $\lambda$. In order to be able to do this, we need to
take care of the following subtle ``uniformity'' issue.
For each local field $E$ individually, we know that
the sum over $\ab$ is finite, by Lemma \ref{fin}. However, we do not know 
whether
this finite set of indices $\ab$ (denoted by $A_{\Gamma}$ in Lemma \ref
{fin}) is {\it the same} for almost all places.
Another difficulty (which would be automatically resolved if we 
knew the affirmative answer to the above question) 
is that for each $\ab$, the formula (\ref{l-dep}) holds
for all but a finite number of places. We need to show that overall 
(i.e. for all $\lambda$ altogether) the set of
places that needs to be eliminated is finite. 
These difficulties are  taken care of by the  
following Theorem due to T.C. Hales \cite[Theorem 2]{toi}. 

First, we need to introduce a notation.
Let $\vartheta(m_1,\dots, m_l)$ be a formula in Pas's language, with
free variables ranging over $\N$.
Let ${\mathcal K}$ be a collection of local fields.
For $E\in {\mathcal K}$, denote by $\vartheta^{E}(m_1,\dots,m_l)$
the interpretation of $\vartheta$
in the model for Pas's language given by the field $E$.

\begin{lem}\label{finie}\cite[Theorem 2]{toi}
Let $\vartheta(m_1,\dots, m_l)$ be a formula in Pas's language.
Suppose that there exists a finite set of prime numbers $S\subset \N$, such
that
for any $E\in {\mathcal K}$, if the residual characteristic of
$E$ is not in $S$, then
$\{(m_1,\dots, m_l)\in\N^l\mid \vartheta^{E}(m_1,\dots,m_l)\}$ is a bounded
subset of $\N^l$.
Then there exists a {\it finite} set of primes $S'$, $S\subset S'\subset\N$, and
exists a {\it bounded} set ${\mathcal C}\subset \N^l$, such that
for any $E\in {\mathcal K}$, if
the residual characteristic of $E$ is not in $S'$,
then
$\{(m_1,\dots, m_l)\in \N^l\mid \vartheta^E(m_1,\dots,m_l)\}\subset
{\mathcal C}$.
\end{lem}

For each element of the finite set of parameters 
corresponding to the conjugacy classes $C$, 
we apply  Lemma \ref{finie} to the formula $\vartheta(\lambda):=$`$\exists x: \phi_{C,\ab}(x)$', 
where $\phi_{C,\ab}(x)$ is the formula of Lemma \ref{conjcl} that defines the subsets
$W_{C,\ab}^E(\Gamma_{\alpha,E})$. That is, $\vartheta^E(\lambda)$
 takes the value `true' iff 
the set $W_{C,\lambda}^E(\Gamma_{\alpha,E})$ 
is not empty (here we understand the symbol $\lambda$
as an abbreviation for the $l$-tuple of variables $(m_1,\dots, m_l)$.)
By Lemma \ref{fin}, for any $E$ such that the notation $W_{C,\lambda}^E$ makes sense 
(that is, for all but finitely many fields $E\in {\mathcal K}$),
the set $\{\lambda\mid\vartheta^E(\lambda)\}$ is finite. Hence, the assumption of 
Lemma \ref{finie} holds. We conclude that 
there exists a finite set of primes $S'$ such that
whenever the residual characteristic of $E$ does not fall in this set,
$\{\lambda\mid\vartheta^E(\lambda)\}\subset {\mathcal C}$, where ${\mathcal C}$ is a bounded
(i.e., finite) subset of $\N^l$.

Finally,  set $M_{\alpha, w}=\sum_{\lambda}M_{\lambda}$, where  
the summation is over the union of finite sets ${\mathcal C}$ 
that correspond to the data defining the conjugacy classes $C$
(which is a finite union of finite sets, so the sum is finite).
This completes the proof.

\section{Appendix: Some Background}
This section is a compilation of brief reviews 
of various concepts and techniques
from logic, algebraic geometry, and representation theory 
that we used above.
\subsection{Pas's language}\label{pl}
We need to deal with families of local fields (e.g. all possible 
nonarchimedean completions of a given number field). Therefore, it is 
necessary to set up a framework that allows to exploit the structure of
a local field without referring to its individual features such as 
the uniformizer of the valuation, for example. 
This is achieved by using a formal language of logic, which has the terms
for the valuation, etc., but does not have the term for a uniformizer.   
This language is called Pas's language. 

A sentence in Pas's
language for a valued field $E$ is allowed to have variables of 
three kinds:
variables running over the valued field $E$ (the {\it valued
field sort}), variables running over its residue field $k$ (the {\it residue field sort}),
and variables running  over $\Z$ (the {\it value sort}).
Formally, it is a three-sorted first order language
\cite{End}.

For variables of the valued field
sort and for variables of the residue field
sort, the language has the
operations of addition (`$+$') and multiplication (`$\times$');
for variables of the value sort (i.e., for $\Z$), only addition
is allowed.
The value sort, additionally, has symbols $\le$, and $\equiv_n$ for
 congruence
modulo each value $n\in\N$.

The language also has symbols for universal (`$\all$') and
existential (`$\exists$') quantifiers, and  standard symbols
$\wedge$, $\vee$, $\neg$, respectively, for logical conjunction,
disjunction, and negation.
We note that the restriction of  Pas's language to the 
residue field sort coincides with the {\it first order language of rings}
\cite{End}.

Naturally, there are symbols denoting all the integers in the value sort.
In the valued field sort, 
there are symbols `$0$' and `$1$', defined as the symbols denoting the
additive and multiplicative  unit,
respectively. Once we have these,
we can formally add to the language the symbols denoting other integers in the
valued field sort. Thus,  `$2$' is the abbreviation for
`$1+1$',  `$-1$' is the abbreviation for `$\exists x, x+1=0$', etc.
Notice that there are no symbols denoting other elements of the valued
field, in particular, there is no symbol for the uniformizer. 
This, of course, agrees with our goal of being able to use the same language 
for any local field within a given family.
We illustrate the allowed operations with the following example.
\begin{ex}
We are not allowed to use any symbols for 
{\it constants} in the valued field.
For example, 
the expression `$a_1x+a_2y=0$' makes sense in Pas's language {\it only}
if $a_1$, $a_2$, $x$ and $y$  are all treated on equal footing, as 
{\it variables} of the same sort. 
\end{ex}

The language also
contains the  following symbols denoting functions:
the symbol `$\val$' for the valuation -- a function from the
valued field sort to the value sort, and 
the symbol `$\ac$' to denote the {\it angular component}
map from the valued field sort to the residue field sort 
(the role of this map will be explained in Section \ref{def}).
We also add the symbol `$\infty$'
to the value sort, to denote the valuation of $0$.

There is a theorem
due to Pas \cite{Pas} that this language  admits quantifier elimination of 
the quantifiers ranging
over the valued field sort and over the value sort
(this is the reason for the absence of multiplication for 
the integers: otherwise, by G\"odel's theorem, there 
would be no quantifier elimination).
This means that every formula in Pas's language can be replaced by an 
equivalent formula without quantifiers ranging over
the valued field sort or over the value sort. 
Further, a theorem due to Presburger \cite{Pres} states that the quantifiers ranging over
the residue field can also be eliminated, so that every formula in Pas's
language can be replaced by an equivalent formula without quantifiers. 
Ultimately, it is this property that makes {\it arithmetic motivic 
integration} possible \cite{DLar}.

\subsection{Definable subsets for $p$-adics}\label{def}
This subsection is, essentially, 
quoted from \cite[Sections 8.2, 8.3]{DLar}.

Let $E$ be the field of fractions of a complete discrete valuation ring
$\ri_{E}$ with finite residue field $k$. 
It is possible to think of $E$ as a {\it structure} 
(in the sense of logic) for Pas's language. That is, we can let the variables
in the formulas in Pas's language  range over $E$ and $k$ respectively,
and then each formula will have true/false value.
In order to match Pas's language with the structure of the field $E$
completely, we need to give a meaning to the symbols that express functions
in Pas's language.
   We  fix a uniformizing parameter $\varpi$. 
The valuation on $E$ is normalized 
so that
$\val(\varpi)=1$.
If $x\in\ri_{E}^{\ast}$ is a unit, there is a natural definition
of $\ac(x)$ -- it is the reduction of $x$ modulo the ideal $(\varpi)$.
Define, for $x\neq 0$ in $E$, $\ac(x)=\ac(\varpi^{-\val(x)}x)$,
 and $\ac(0)=|0|=0$.

Now let us take $F$ -- a finite extension of $\Q$ with 
ring of integers
${\mathcal O}_F$. Following \cite{DLar}, we 
let $R=\ri[\frac1N]$, for some non zero
integer $N$ which is a multiple of the discriminant of $F$.
Let $v$ be a closed point of $\spec R$. We denote by $F_v$ the 
completion of
the localization of $R$ at $v$, by $\ri_{F_v}$ -- its ring of integers,
and by $k_v$ -- the residue field. If $N$ equals 
the discriminant of $F$, 
then $\{F_v\}$ is the family of all completions of $F$ at 
the places lying over  the primes in $\Q$ that do not ramify
in $F$.

In Section \ref{pl}, we formally added to Pas's language 
the symbols to denote all 
integers in the valued field sort. In order to work with the family 
of fields ${\mathcal K}=\{F_v\}$, it is convenient to also add
to Pas's language, for every element of the global field $F$, 
a symbol to denote this element in the valued
field sort. This allows to ``have constants from $F$'' in the formulas
in Pas's language. 
We will call this extended language the Pas's language for the
valued field $F((t))$. The reason for this is that $F((t))$ is, naturally, 
a valued field ($t$ being the uniformizer of the valuation), and every logical
formula in the extended Pas's language for this field also makes sense as 
a formula in the extended Pas's language for every completion of the field $F$.
It is through interpreting the Pas's language formulas in the model given by 
the field $F((t))$ that one eventually arrives at geometric objects 
associated with them (see \cite[Section 6.1]{TH1}). 

\begin{defn}
Let $\phi$ be a formula in the  Pas's language for the valued field
$F((t))$, with $m$ free variables running over the valued field sort and
no free variables running over the residue field sort or the value sort.
For each $v$ -- a place of $F$, and $E=F_v$ -- the corresponding completion,
denote by $Z(\phi,\ri_{E})$
the subset of $\ri_{E}^m$ defined by the formula $\phi$:
$Z(\phi, \ri_E)=\{(x_1,\dots,x_m)\in\ri^m\mid \phi(x_1,\dots, x_m)\}$
(recall that $\phi$ has the values true/false). 
A subset $B\subset \ri_{E}^m$ is called {\it definable} if there exists
a formula $\phi$ such that $B=Z(\phi, \ri_{E})$.
\end{defn}

For example, if $V\hookrightarrow \A^m$ is an affine variety,
$V(\Z_p)$ is a definable subset of $\Z_p^m$ (it can be defined by a
 logical formula with $m$ free variables which does not use the symbols 
`$\ac$' and `$\val$': just the polynomial relations defining the variety 
$V$). 

So far, we have defined the notion of a definable subset of
$\ri_{E}^m=\A^m(\ri_{E})$. 
This definition extends naturally to give a notion of a {\it definable}
subset of an affine variety $\mx$ over $\ri_{E}$.

Let $\mx$ be a smooth variety over $F$ of dimension $d$. 
There is a natural
$d$-dimensional measure on $\mx(\ri_{F_v})$ for each $v$, 
which we shall denote by
$vol_{\mx,F_v}$. 
This is the Serre-Oesterl\'e measure, \cite{oe}, \cite{se}.
All definable subsets of $\mx(\ri_{F_v})$ 
are measurable with respect to this
measure, when $F$ has characteristic $0$, \cite{den}. This measure is
defined by requiring that the fibers of the reduction map modulo
$\varpi_v$ have volume $q^{-d}$, where $q$ is the cardinality of 
the residue field $\ri_{F_v}/(\varpi_v)$.

In the Section \ref{armint}, we describe the theory of 
{\it arithmetic motivic integration} which associates 
algebraic geometric objects with logical formulas.
For each place  $v$, 
there is a linear operator (Frobenius) acting on the cohomology
of these objects, and its trace on the object associated with the formula
$\phi$ equals the Serre-Oesterl\'e volume of the set $Z(\phi, \ri_{F_v})$
for almost all $v$. 

\subsection{Motivic measures}\label{asp}
The purpose of this subsection is twofold.
The first goal is to outline the context for the ideas 
from the theory of motivic
integration that we are using and to give some background references.
The second objective is to let someone who is an  expert in motivic 
integration know exactly to what extent the theory is being used in this 
paper.
All the precise definitions and statements that we need are  
quoted in the next three subsections.

The term  {\it motivic integration} first appeared in M. Kontsevich's 
lecture at Orsay in 1995. Let $k$ be an algebraically closed field of 
characteristic  $0$, and let $\mx$ be an algebraic variety over $k$.
{\it Motivic measure} lives on the {\it arc space} of $\mx$, i.e. on the 
``infinite-dimensional variety'' whose set of $k$-points coincides with the
set of $k[[t]]$-points of $\mx$. 
It is a measure in every sense 
(i.e., it is additive, it transforms under 
morphisms in a such a way that makes it analogous to measures on real or
$p$-adic manifolds, etc.),
but its nature is algebraic, and its values are not real or complex numbers,
but, roughly speaking, equivalence classes of varieties.
The theory of motivic integration is, in fact, a theory about the motivic 
measure. The functions [on the arc space] 
that can be integrated against this measure are 
scarce.
Due to algebraic nature of the measure, the algebra of measurable subsets
of the arc space is rather coarse: it is, essentially, generated by
the sets ``growing out'' of subvarieties of $\mx$.
For  background information on the original motivic integration we refer
to \cite{craw}, \cite{DL}, \cite{DL1},  and  
to the Bourbaki talk by E. Looijenga \cite{Loj}.

The {\it arithmetic motivic integration} developed 
by J. Denef and F. Loeser in 1999 \cite{DLar}
uses similar ideas, but it is an independent theory. 
The main difference is that it is adapted to deal with sets of rational 
points of the given variety over various fields, as opposed to the 
less flexible original theory which is suited only to deal with the 
arc space of the given variety, and not its rational points.  
The {\it arithmetic motivic volume} 
takes values in {\it virtual Chow motives}.

Instead of the arc space itself, the arithmetic motivic volume
lives on 
{\it definable subassignments} of the functor of its points.
This distinction allows one to work directly at the level of rational points 
\cite{DLar}.
The {\it definable subassignments}
are, in a sense, a geometric incarnation of 
logical formulas in Pas's language.
We will not need these geometric objects, since all the varieties we 
deal with here are smooth and affine. Hence, we will assign ``motivic 
volumes'' to logical formulas directly.

\subsection{The ring $\mot$}\label{mots}
Here we define the ring of values
of the arithmetic motivic volume. The {\it volume} itself will be defined
in the next subsection. 
This section is, essentially, quoted from \cite[Section 1.3]{DLar}, with one modification: 
the definition of the ring of values was greatly simplified by Denef and Loeser  in 
\cite{DLn}. Hence, we are using the version of this ring that appears in \cite{DLn}
rather than the original one.

Let $F$ be a field.
Let us denote by ${\text {Mot}}_F$ the category of Chow motives over $F$,
with coefficients in $\Q$.
The objects in ${\text {Mot}}_F$ are, formally, equivalence classes of
triples $(S,p,n)$,
where $S$ is a proper and smooth scheme over $F$, $p$ is an idempotent
correspondence on $S$ with coefficients on $\Q$
[that is, an algebraic cycle in $S\times S$ such that $p^2\simeq p$, 
where $p^2$ is the product of $p$ with itself], and $n\in \Z$.
Some details about Chow motives can be found in \cite{scholl}.

Even though Chow motives do not form an abelian category, they are 
sufficiently suited for our purposes, because the category 
${\text {Mot}}_F$ is pseudo-abelian, and therefore the notion 
of its Grothendieck group makes sense. 
We denote this Grothendieck
group by $K_0({\text {Mot}}_F)$. This is the additive group of 
equivalence classes of formal linear
combinations of objects of ${\text {Mot}}_F$ (with the natural notion of 
equivalence so that if $M=A\oplus B$ then $[M]=[A]+[B]$). 
Since the category of Chow motives has
a tensor product, $K_0({\text {Mot}}_F)$ can be given 
the structure of a ring.
There is a natural inclusion of the set of objects of the category 
$\text{Mot}_F$ into its Grothendieck ring (each object can be identified with 
its own
equivalence class in  $K_0({\text {Mot}}_F)$).

Let ${\text {Var}}_F$ be the category of algebraic varieties over $F$.
For a {\it smooth projective} variety $S$, there is a Chow motive canonically
associated
with it  -- the triple $(S, id, 0)$.
There exists
a unique embedding 
$$
\chi_c:{\text {Var}}_F\to {\text {Mot}}_F
$$
such that for smooth projective $S$, $\chi_c(S)=(S,id,0)$.\newline
Let $\lef=(\spec F, id, -1)$ be the Lefschetz motive
(see \cite{scholl} for details).
It is the image, under the morphism $\chi_c$, of  
the affine line
$\A^1$ (note that the affine line is not a projective variety, and 
this is the reason the integer component of the Chow 
motive associated with it is not $0$).
The map $\chi_c$ induces the morphism of Grothendieck rings 
\cite{GS}, \cite{NA}: 
$$
\chi_c:K_0({\text {Var}}_F)\to K_0({\text {Mot}}_F).
$$
The object $\lef$ is invertible (with respect to tensor product) in 
the category of Chow motives.
Hence, the map $\chi_c$ extends to the localization
of $K_0({\text {Var}}_F)$ at $[\A^1]$.
We denote by $K_0^{mot}(\text{Var}_F)$ the image of 
this extended map $\chi_c$.
Finally, set $K_0^{mot}({\text {Var}}_F)_{\Q}=K_0^{mot}(\text{Var}_F)\otimes \Q$.

For more details of the  definition of the ring of values of 
the arithmetic motivic volume we refer to
\cite[Section 1.3]{DLar} and \cite[Section 6.3]{DLn}.
In the next subsection, we describe the 
{\it arithmetic motivic volume}, which
is a function on logical formulas [rather, on {\it definable 
subassignments} as
in \cite{DLar}] taking values, roughly, in this ring.
More precisely, it is necessary to complete the 
ring $K_0^{mot}({\text {Var}}_F)_{\Q}$ with respect to 
a certain
{\it dimensional filtration} in order to have a 
meaningful ``measure theory''
with values in it.  To define the filtration of 
the ring of virtual Chow motives, we first define  a filtration on 
$K_0({\text {Var}}_F)_{loc}$ (where the subscript `$loc$' stands for 
localization at $[\A^1]$). 
For $m\in \Z$, let 
$F^m K_0({\text {Var}}_F)_{loc}$ 
be the subgroup of $K_0({\text {Var}}_F)$ 
generated by the elements of the form $[S]\lef^{-i}$ with $i-\dim S\ge m$.
This defines a decreasing filtration $F^m$. The image of this filtration 
under the map $\chi_c$ is a decreasing filtration on 
$K_0^{mot}(\text{Var}_F)$, which naturally induces a filtration on the full ring
$K_0^{mot}({\text {Var}}_F)_{\Q}$.
The completion of this ring with respect to this 
filtration is denoted by $\mot$. 

\begin{rem}It is important to note that 
the volumes of all the objects that we will be
dealing with are contained in the image of the ring 
$K_0^{mot}({\text {Var}}_F)_{\Q}$
in $\mot$  under the  completion map, if we adjoin to it all the 
elements of the 
form  $[\lef^i-1]^{-1}$ for positive integers $i$ 
(see \cite[Remark 8.1.2]{DLar}).
In fact,  there is now new, yet unpublished, version of the theory 
of motivic integration (due to  R. Cluckers and F. Loeser), 
which does not use the completion at all. 
This new theory allows a refinement of our results, namely,
the virtual Chow motives responsible for the character values should 
automatically lie in the ring 
$K_0^{mot}({\text {Var}}_F)_{\Q}[\lef^i-1]_{i\ge 1}$, 
rather than in a less understood complete ring $\mot$. 
\end{rem}

If $v$ is a place of $F$, there is an action of $Frob_v$ on the 
elements of $\mot$ that comes from the Frobenius action on the 
Chow motives. The trace of the
Frobenius operator on a Chow motive is the alternating sum 
of the traces of Frobenius acting on its $l$-adic cohomology 
groups, and it is an element of $\bar {\Q_l}$ (where $l$ is a prime 
number), see \cite[Section 3.3]{DLar}.  
In all the cases that we consider (following Denef and Loeser),
this number turns out to lie in  $\Q$   
(in particular, the choice of $l$ doesn't matter). It is the trace of 
Frobenius action that allows to relate the values of the motivic 
measure, that are elements of $\mot$, with the usual $p$-adic volumes, 
that are rational numbers.

\subsection{Arithmetic motivic integration}\label{armint}
A beautiful very short description of the theory of 
arithmetic motivic integration can be found in \cite[Section 6.1]{TH1}.

Let $\mx$ be a smooth $d$-dimensional affine 
variety over $F$, where $F$ is a global field
as above. Let $\phi$ be a logical formula in Pas's language for the
field $F((t))$ with $m$ free variables ranging over the valued field
sort and no bound variables over the value sort.
There is a natural notion of the formula $\phi$ defining a subset of 
$\mx$: suppose that $\mx$ is embedded in $\A^m$ in some way. Then 
we can assume that $\phi(x_1,\dots, x_m)$ is true only if $x_1,\dots,x_m$
satisfy the polynomial equations defining $\mx$. 
In other words, this happens when $Z(\phi, \ri_{F_v})\subset \mx(F_v)$  
for almost all places $v$ of $F$, see 
\cite[Proposition 5.2.1 and Section 8.3]{DLar}. 

In \cite[Section 6.3]{DLn} (following the original method  of \cite[Section 3]{DLar}), Denef and Loeser 
associate elements of the 
ring $K_0^{mot}(Var_F)_{\Q}$  to such formulas $\phi$. This mapping 
is additive in a natural sense (see \cite[Proposition 3.4.4]{DLar}).
 This map is unique; it is denoted by $\chi_c(\phi)$ in \cite{DLar}, but we will
denote it by $\mu(\phi)$. This map eventually extends to the 
formulas that involve quantifiers. In order to define this extension,
the authors use the language of  
definable subassignments of the functor of points of 
$\arcs(\A^m)$. The resulting map on subassignments is 
called {\it arithmetic motivic volume},
see  \cite[Section 6]{DLar}.
We will not need the precise definition of the arithmetic motivic volume 
or the language of subassignments. We only need its existence, 
finite additivity,
and the following two theorems. In the form that we state these theorems,
they are special cases of Theorems 8.3.1 and 8.3.2 \cite{DLar}, 
respectively, corresponding to  the case when the variety $\mx$ is smooth 
and affine, and the definable 
subassignment is defined by just one formula on the whole space, 
in the terms of \cite{DLar}. 

\begin{thm}\cite[Theorem 8.3.1]{DLar}
Let $F$ be a finite extension of $\Q$ with the ring of integers $\ri$,
and $R=\ri[\frac1N]$, for some non zero integer $N$.
Let $\mx$ be  an affine  variety over $R$, and let $\phi$ 
be a logical formula as above. 
Then there exists a non zero multiple $N'$ of $N$,
such that, for every closed point $v$ of $\spec \ri[\frac1{N'}]$,
$$
\tf_v\bigl(\mu(\phi)\bigr)=vol_{\mx,F_v}(Z(\phi,\ri_{F_v})).
$$
\end{thm}

In the case when the local fields under consideration
are of finite characteristic, 
there is a slightly stronger result:

\begin{thm}\cite[Theorem 8.3.2]{DLar}
Let $F$ be a field of characteristic $0$ which is the field of fractions of a
normal domain $R$ of finite type over $\Z$. Let ${\mathcal X}$ be a
variety over $R$, and
let $\phi$ be a logical formula.
Then there exists a nonzero element $f$ of $R$, such that, for every closed point
$x$ of $\spec R_f$, $Z(\phi,\F_x[[t]])$ is $vol_{\F_x[[t]]}$-measurable
and
$$
\tf_x\bigl(\mu(\phi)\bigr)=vol_{\mx,F_x[[t]]}(Z(\phi,{\F_x[[t]]})).
$$
\end{thm}

This completes our survey of arithmetic motivic integration. 
The next two subsections are devoted to a summary of the facts from 
representation theory that we are using.

\subsection{Deligne-Lusztig representations of classical groups}\label{delu}
Let $\G$ be a connected reductive algebraic group defined over
$\bar \F$, where $\F=\F_q$ is a finite field, $q=p^m$.
Deligne and Lusztig \cite{DeLu} constructed a large class of representations
of groups of the form $\G^{\F}$ in vector spaces
 over $\bar{\Q_l}$ with $l\neq p$.
By $\G^{\F}$ we denote the group of fixed points of $\G(\bar\F)$
under the
Frobenius map associated with $q$.

These representations are parametrized by
pairs $(T, \chi)$ where $T$ is a maximal torus in $\G^{\F}$ 
defined over
$\F$ and $\chi$ is an irreducible character of $T$. A representation
corresponding to the pair $(T,\chi)$ is denoted by
$R_{T,\chi}^{\G^{\F}}$.
\begin{defn}
An irreducible representation of $\G^{\F}$ is called {\it Deligne-Lusztig}
if it is equivalent to $R_{T,\chi}^{\G^{\F}}$ for some $(T,\chi)$.
\end{defn}

In this paper we  deal with representations of
$p$-adic groups that are obtained  from Deligne-Lusztig representations
by a certain natural procedure described in Section \ref{inflation}.
The reason we restrict our attention only
to Deligne-Lusztig representations is the fact that the values of their
characters at unipotent elements in ${\G}^{\F}$ are given by {polynomials} in
$q$. We use this fact in an essential way in the proof of 
our main result.

More specifically, if $u$ is a unipotent element of $\G^{\F}$,
the value of the character of $R_{T,\chi}^{\G^{\F}}$ at 
$u$ can be expressed as a
sum of values of Green functions $Q_T^{\G}(u)$ (see, e.g.,
\cite[Theorem 6.8]{Sri}) (in particular, it does not depend on $\chi$).
It was proved by B. Srinivasan in \cite{Sr3} that
the values of $Q_T^{\G}(u)$ are polynomials in $q$ (the polynomial
itself depends on the conjugacy class of $u$) 
if $\G$ is symplectic or odd special orthogonal,
and the characteristic $p$ is odd (since everything we do here is only
up to a finite number of primes, this is not a restriction).

\subsection{Parametrization of unipotent conjugacy classes}\label{unip}
Here we review the parametrization of nilpotent orbits in the classical
$p$-adic
Lie algebras given by Waldspurger \cite{W}.
We quote it from \cite[Sections I.5, I.6]{W}.
It is used in an essential way in the proof of Lemma \ref{conjcl}.

Let $\lie$ be the Lie algebra $\syp(r)$ or $\so(r)$ with $r$ odd, and
let $X\in \lie$ be a nilpotent element.
Denote by $(V,q_V)$ the underlying vector space of $\lie$ with the
quadratic form from the definition of  $\lie$ on it.

Consider the set of all partitions $\Lambda=(\Lambda_j)$
of $r$ with the following properties:
\begin{itemize}
\item
In the symplectic case, for all odd $i\ge 1$, $c_i(\Lambda)$ is even
\item
In the orthogonal case, for any even $i\ge 2$, $c_i(\Lambda)$ is even
\end{itemize}
(here $c_i(\Lambda)$ is the number of $\Lambda_j$'s that  equal $i$).

We can associate with $X$ a partition $\Lambda$ of $r$: for all integers
$i\ge 1$, $c_i(\Lambda)$ is the number of Jordan blocks of $X$ of length $i$
in the natural matrix representation.
This partition automatically  satisfies the above conditions.

For all $i\ge 1$, set
$V_i=\ker(X^i)/[\ker (X^{i-1})+X\ker(X^{i+1})]$.
In the symplectic (resp, orthogonal) case, define the quadratic form
$\tilde q_i$ on $\ker(X^i)$, for all even $i$ (resp., odd), by:
$$
\tilde q_i(v,v')=(-1)^{\left[\frac{i-1}2\right]}q_V(X^{i-1}(v),v').
$$
Passing to a quotient, we get a non-degenerate form $q_i$ on $V_i$.

In the orthogonal case, the forms $q_i$ satisfy the condition
$$
\bigoplus_{i {\  \text{odd}}} q_i\sim_a q_V,
$$
where $\sim_a$ indicates that the two forms have 
the same anisotropic kernel.

The correspondence described above is a bijection between the set of
nilpotent orbits (under conjugation by elements of $G$) in $\lie$ and
the set of pairs $(\Lambda, (q_i))$ with $\Lambda$ -- a partition of $r$
and $(q_i)$ -- a collection of quadratic forms of dimensions
$c_i(\Lambda)$ satisfying the conditions mentioned above.

The same classification holds for finite fields.
We also observe that over a finite field $\F_q$,
an equivalence class of a quadratic
form is determined by its rank and discriminant in
$\F_q^{\ast}/{\F_q^{\ast}}^2$
\cite{Serre}.
The nilpotent orbits in the Lie algebra are in bijection with the unipotent
conjugacy classes in the group, when the characteristic $p$ is large enough
\cite{hump}. Explicitly, the bijection is given by  Cayley transform
$x\mapsto(1-x)(1+x)^{-1}$.

Let $\G$ be a simply connected
algebraic group of symplectic or odd orthogonal type.

Finally, we see that the unipotent conjugacy classes of $\G^{\F_q}$
are in bijection with the following data (cf. also \cite[I 2.9]{Sp} or \cite[11.1]{Lu}):
\begin{itemize}
\item a partition $\Lambda$ of $r$ with certain properties;
\item a representative (one of the two possible choices) of
$\F_q^{\ast}/{\F_q^{\ast}}^2$ for each even (resp. odd) $c_i(\Lambda)$.
\end{itemize}

%%%%%%%%%%%%%%%%%%%%%%%%%%%%%%%%%%%%%%%%
\end{document}